\documentclass{amsart}     
\usepackage{amssymb} 
\usepackage{amsmath,amscd,eucal,latexsym} 
\usepackage{graphicx}
\usepackage{epsf}

\newtheorem{theorem}{Theorem}
\newtheorem{lemma}[theorem]{Lemma}

\newtheorem{prop}[theorem]{Proposition}

\theoremstyle{definition} 
\newtheorem*{rem}{Remark}

\newtheorem{example}{Example}

\DeclareMathOperator{\intr}{int}

\newcommand{\bd}{\partial}

\newcommand{\sm}{\smallsetminus}

\newcommand{\TT}{\mathcal{T}} 
 
\newcommand{\UU}{\mathcal{U}}

\newcommand{\Z}{\mathbb{Z}}

\begin{document}

\title[Endperiodic Automorphisms]{Examples of Endperiodic Automorphisms}

\author[J. Cantwell]{John Cantwell}
\address{Department of Mathematics\\ St. Louis University\\ St. 
Louis, MO 
63103}
\email{CANTWELLJC@SLU.EDU}

\author[L. Conlon]{Lawrence Conlon}
\address{Department of Mathematics\\ Washington University, St. 
Louis, MO 
63130}
\email{LC@MATH.WUSTL.EDU}


\begin{abstract}
We   give examples of endperiodic automorphisms.
\end{abstract}

\maketitle

These examples are intended to supplement the examples in our paper~\cite{cc:HM}. Also see Fenley~\cite{fe:thesis} and~\cite[Section~5]{fe:endp} for other examples. We use the terminology and notation of~\cite{cc:HM}. In the figures, traintracks for the positive laminations $\Lambda_{+}$ are giving with solid curves and traintracks for the negative laminations $\Lambda_{-}$ are given with dashed curves.

\section{Examples with $\Lambda_{\pm}$ finite}

Examples~\ref{ex2leaves} and~\ref{ex3leaves} have the properties that,  

\begin{enumerate}

\item Both $\Lambda_{+}$ and $\Lambda_{-}$ are finite\upn{;}

\item $\Lambda_{+}$ (respectively $\Lambda_{-}$) has a leaf $\lambda$ such that there is an arc $\alpha$ tranverse to $\lambda$ which meets $\Lambda_{+}$ (respectively $\Lambda_{-}$) in a countably infinite set\upn{;}\label{item2}

\item Dehn twists are not used in the construction of the examples.\label{item3}

\item There do not exist transverse measures to $\Lambda_{+}$ and $\Lambda_{-}$with full support.\

\end{enumerate}
Properties~\ref{item2} and~\ref{item3} are surprising. These two examples give non-trivial examples in which both $\Lambda_{+}$ and $\Lambda_{-}$ are finite. We know of no examples in which $\Lambda_{+}$ or $\Lambda_{-}$ are countably infinite. Fenley~\cite[Section~5]{fe:endp} also gives an example of laminations which do not support a measure of full support.

Example~\ref{exbasic} is basic to the construction of Examples~\ref{ex2leaves} and~\ref{ex3leaves}.

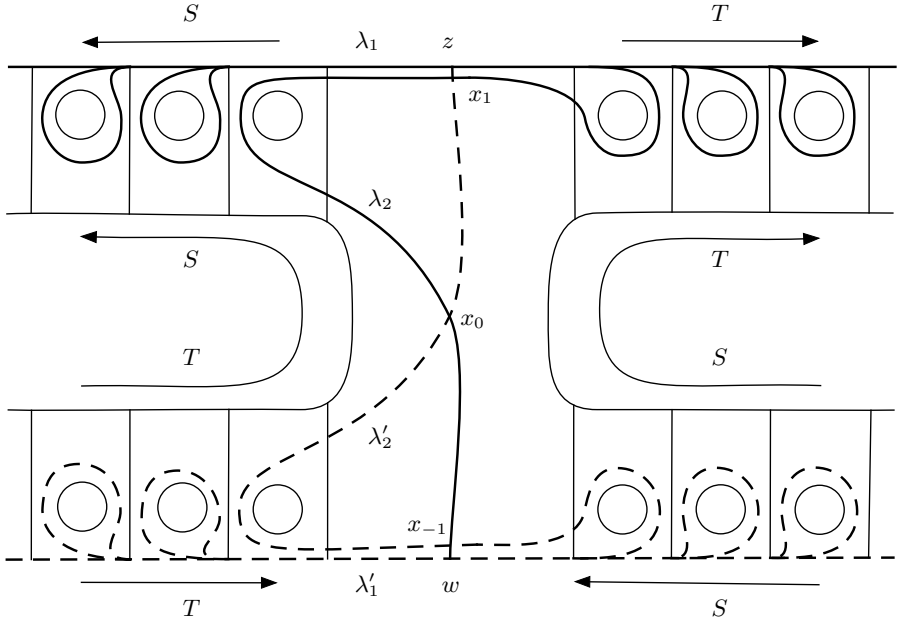
\begin{figure}[h]
\begin{center}
\begin{picture}(300,240)(20,-220)
\rotatebox{270}{\scalebox{.7}{\includegraphics[width=300pt]{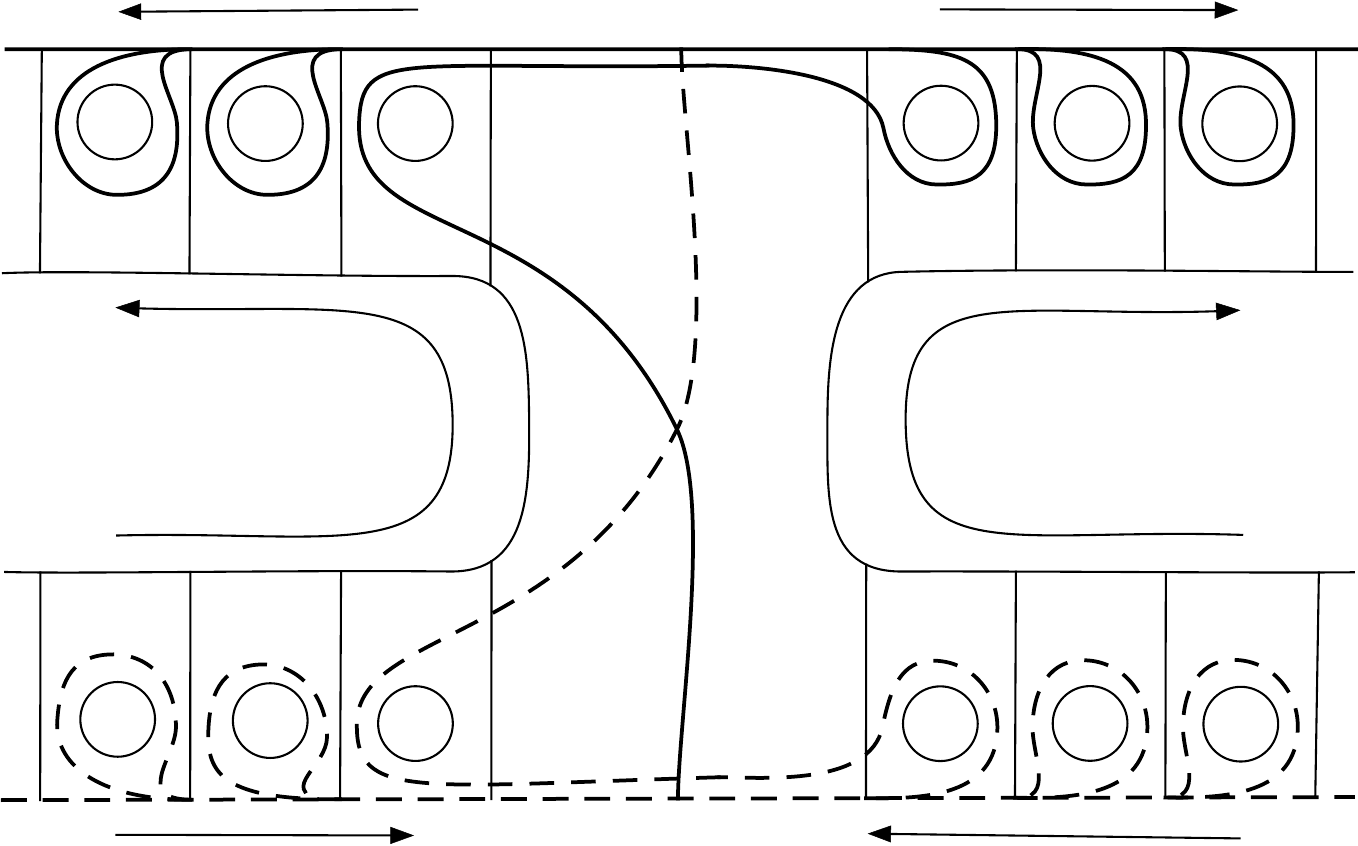}}}
\put (-270,5){\small$S$}
\put (-70,5){\small$T$}
\put (-270,-88){\small$S$}
\put (-70,-88){\small$T$}
\put (-200,-65){\small$\lambda_{2}$}
\put (-205,-5){\small$\lambda_{1}$}
\put (-200,-155){\small$\lambda'_{2}$}
\put (-205,-211){\small$\lambda'_{1}$}
\put (-185,-188){\small$x_{-1}$}
\put (-165,-110){\small$x_{0}$}
\put (-162,-25){\small$x_{1}$}
\put (-172,-5){\small$z$}
\put (-172,-211){\small$w$}
\put (-270,-125){\small$T$}
\put (-70,-125){\small$S$}
\put (-270,-220){\small$T$}
\put (-70,-220){\small$S$}
\end{picture}
\caption{The surface $L$ and traintracks for $\Lambda_{\pm}$ in Example~\ref{exbasic}}\label{fig1}
\end{center}
\end{figure}

\begin{example}\label{exbasic}

Let $L$ be  the  surface in Figure~\ref{fig1}.    On $L$, the homeomorphism $T$ moves the circles in the bottom left leg of  $L$ (and in fact the fundamental domains) to the right  by one, moves the rightmost circle in the bottom left leg to the leftmost circle in the top right leg of $L$,  moves the circles (and fundamental domains) in the top right leg of $L$ to the right by one, and is pointwise fixed on the top left and bottom right legs of $L$.  Similarly, the homeomorphism $S$ moves the circles (and fundamental domains) in the bottom right leg  of $L$ to the left by one, moves the leftmost circle in the bottom right leg of $L$ to the rightmost circle in the top left leg of $L$,  moves the circles (and fundamental domains) in the top left leg of $L$ to the left by one, and is pointwise fixed on the top right and bottom left legs of $L$.  The endperiodic automorphism $f = S\circ T:L\to L$ behaves as indicated on the four real line boundary components of $L$. In particular, $f$ has a fixed point $z$ on the top  edge of $L$  and a fixed point $w$ on the bottom  edge of $L$.

The lamination $\Lambda_{+}$ has two leaves, $\lambda_{1}$ is the top edge of $L$ in Figure~\ref{fig1} and has two escaping ends. The other leaf $\lambda_{2}$ is a semi-isolated half-line that approaches $\lambda_{1}$ from below. Similarly, the lamination $\Lambda_{-}$ has two leaves, $\lambda'_{1}$ is the bottom edge of $L$ in Figure~\ref{fig1} and has two escaping ends. The other leaf $\lambda'_{2}$ is a semi-isolated half-line that approaches $\lambda'_{1}$ from above. The intersection $\lambda_{2}\cap\lambda'_{2}$ consists of points $x_{n}$, $n\in\Z$ which can be indexed so that,

\begin{enumerate}

\item The sequence $\{x_{n}\}$ is monotone  in both $\lambda_{2}$ and $\lambda'_{2}$\upn{;}

\item $x_{n}\to z$ as $n\to+\infty$ and  $x_{n}\to w$ as $n\to-\infty$\upn{;}

\item $f(x_{n}) = x_{n+1}$, $n\in\Z$.

\end{enumerate}
Note $x_{-1},x_{0},x_{1}$ in Figure~\ref{fig1}. The set $\Lambda_{+}\cap\Lambda_{-}$ consists of the points $z$, $w$, and $x_{n}$, $n\in\Z$.

\end{example}

\begin{figure}[h]
\begin{center}
\begin{picture}(300,240)(20,-220)
\rotatebox{270}{\scalebox{.7}{\includegraphics[width=300pt]{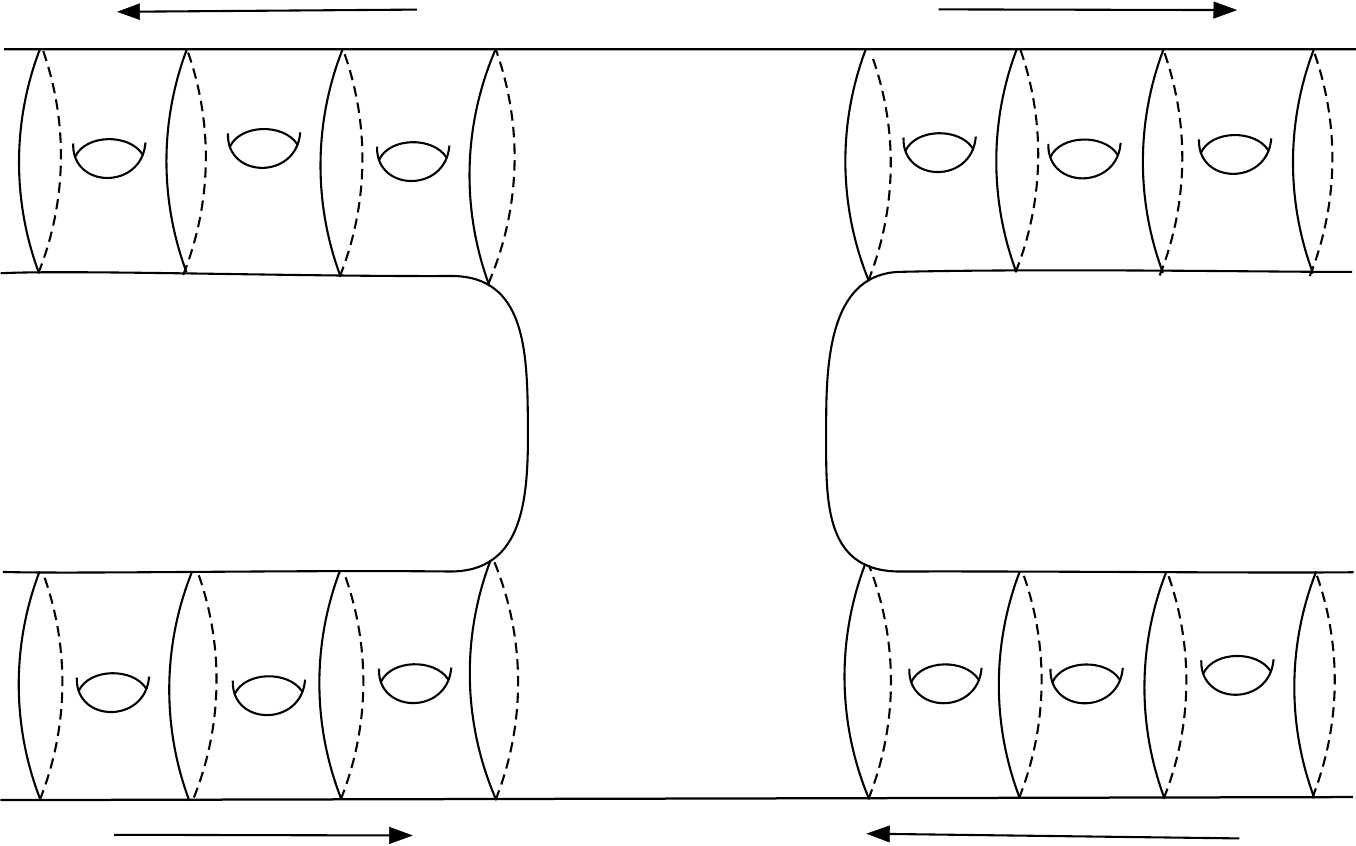}}}
\put (-270,5){\small$S$}
\put (-70,5){\small$T$}
\put (-270,-220){\small$T$}
\put (-70,-220){\small$S$}
\end{picture}
\caption{The surface $L$ in Example~\ref{ex2leaves}}\label{four1}
\end{center}
\end{figure}
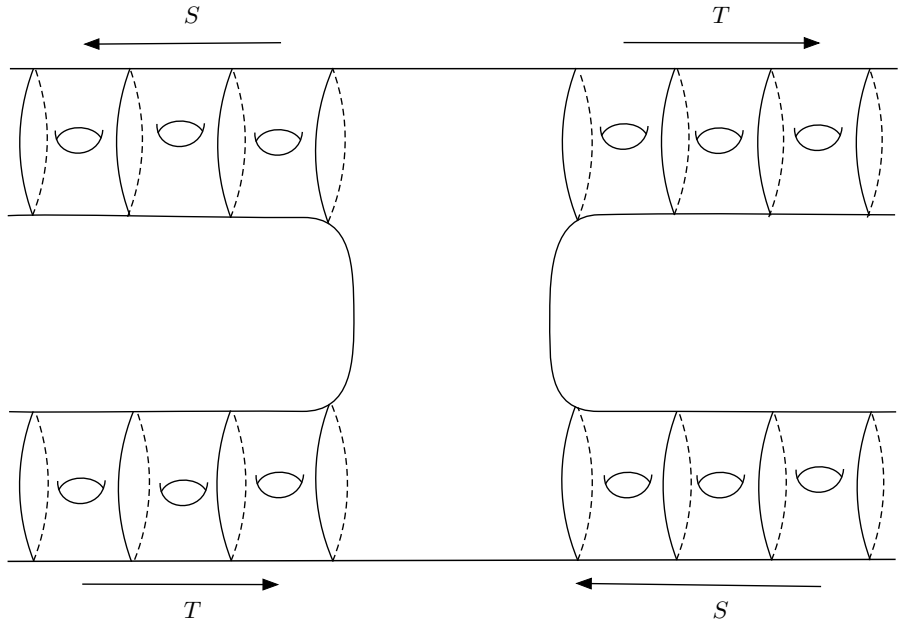

\begin{example}\label{ex2leaves}

Let $L$ be the double of the suface of Example~\ref{exbasic} and $f$ the double of the endperiodic map. The surface $L$ is depicted in Figure~\ref{four1} and has four non-planar ends. The lamainations $\Lambda_{\pm}$ will be the doubles of the laminations in Example~\ref{exbasic}. The traintracks on Figure~\ref{four1} can be visualized from the traintracks in Figure~\ref{fig1}.  In particularl, $\Lambda_{+}$ has two leaves. The first $\lambda_{1}$ has two escaping ends and is approached on each  side by an end of the second $\lambda_{2}$.  The leaf $\lambda_{2}$ is isolated and returns infinitely often to the core. In Figure~\ref{four1} the leaf $\lambda_{1}$ runs along the top of the figure and one end of $\lambda_{2}$ approaches $\lambda_{1}$ on the front of the figure and the other end of $\lambda_{2}$ approaches $\lambda_{1}$ along the back of the figure. A similarly description can be given of the two leaves $\lambda'_{1},\lambda'_{2}$ of $\Lambda_{-}$. The leaf $\lambda'_{2}$ meets $\lambda_{1}$ in a point $z$. The leaf $\lambda_{2}$ meets $\lambda'_{1}$ in a point $w$.   The intersection $\lambda_{2}\cap\lambda'_{2}$ consists of two   sequences $x_{n}$, $n\in\Z$, and $x'_{n}$, $n\in\Z$, the first on the front of the figure and the second on the back. The indexing can be chosen so that 

\begin{enumerate}

\item The sequences $\{x_{n}\}$ and $\{x'_{n}\}$ are monotone  in both $\lambda_{2}$ and $\lambda'_{2}$\upn{;}

\item $x_{n}\to z$, $x'_{n}\to z$ as $n\to+\infty$ and $x_{n}\to w$, $x'_{n}\to w$ as $n\to-\infty$\upn{;}

\item $f(x_{n}) = x_{n+1}$ and $f(x'_{n}) = x'_{n+1}$, $n\in\Z$.

\end{enumerate}
The set $\Lambda_{+}\cap\Lambda_{-}$ consists of the points $z$, $w$, and $x_{n}$ and $x'_{n}$, $n\in\Z$.

\end{example}

If the leaf $\lambda_{2}$ has positive transverse measure, then any arc transverse to the leaf $\lambda_{1}$ would have infinite measure. Thus there does not exist a transverse measure on $\Lambda_{+}$ of full support.

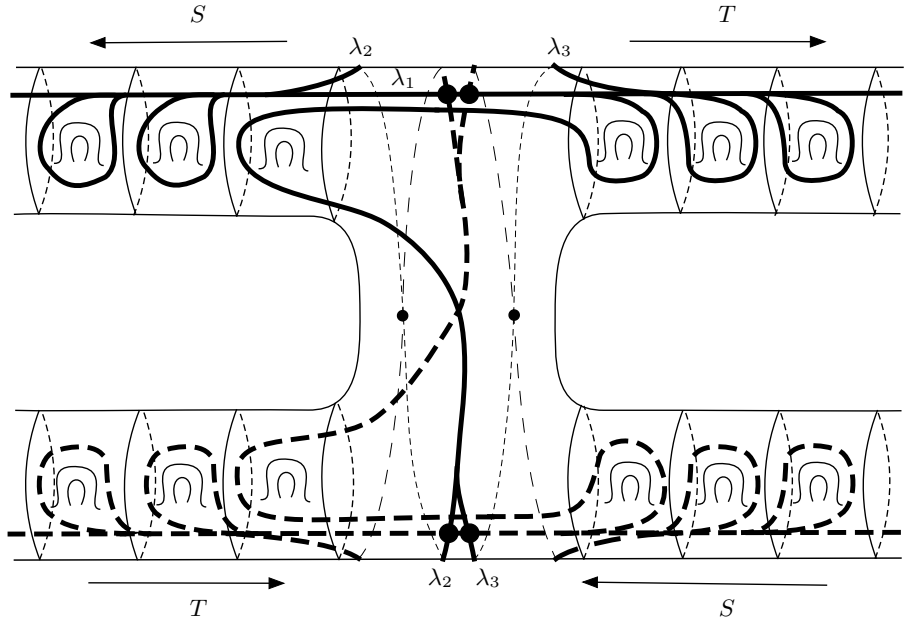
\begin{figure}[t]
\begin{center}
\begin{picture}(300,240)(20,-220)
\rotatebox{270}{\scalebox{.7}{\includegraphics[width=300pt]{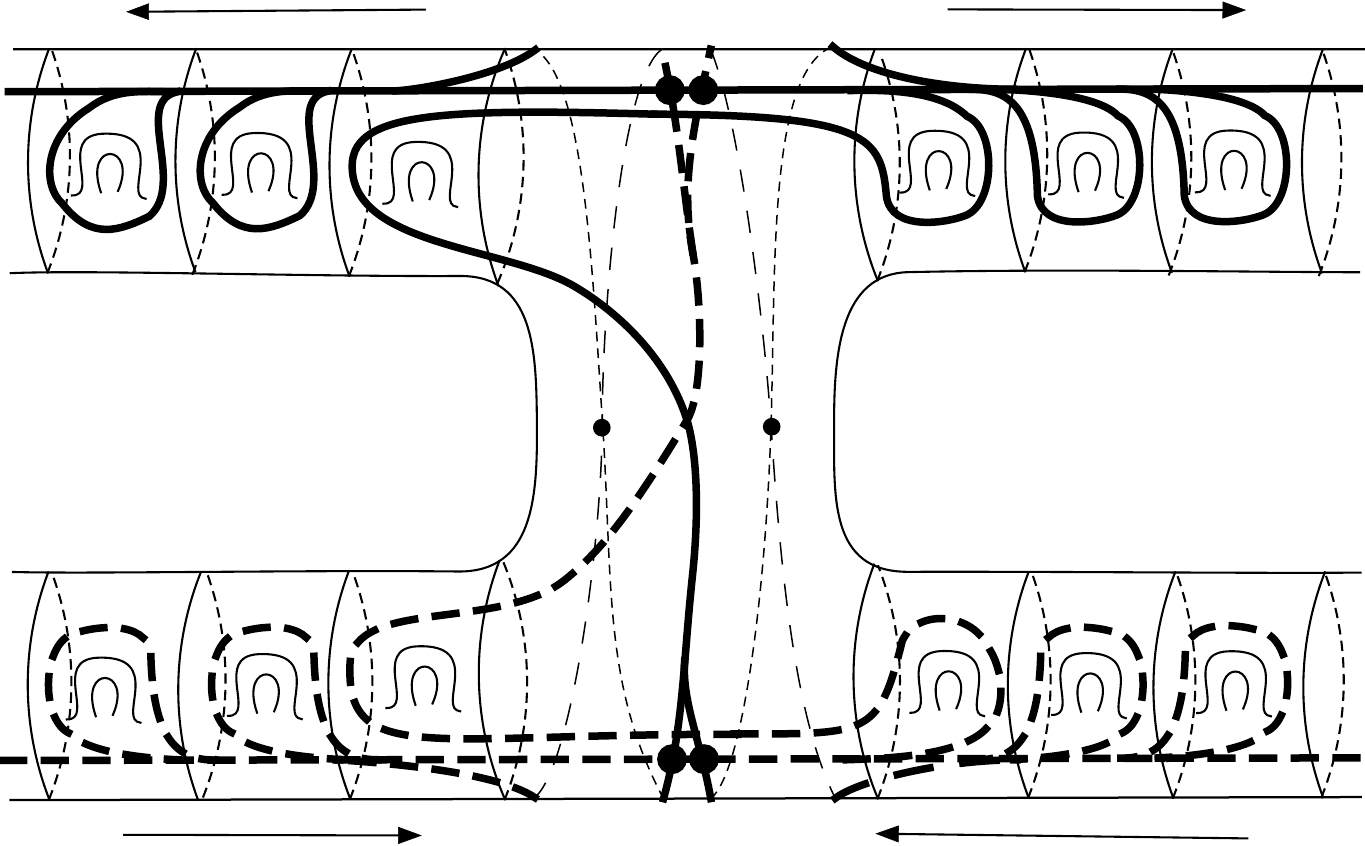}}}
\put (-270,5){\small$S$}
\put (-70,5){\small$T$}
\put (-194,-19){\small$\lambda_{1}$}
\put (-210,-7){\small$\lambda_{2}$}
\put (-135,-7){\small$\lambda_{3}$}
\put (-180,-208){\small$\lambda_{2}$}
\put (-162,-208){\small$\lambda_{3}$}
\put (-270,-220){\small$T$}
\put (-70,-220){\small$S$}
\end{picture}
\caption{The surface $L$ and traintracks for $\Lambda_{\pm}$ in Example~\ref{ex3leaves}}\label{four2}
\end{center}
\end{figure}

\begin{example}\label{ex3leaves}

Example~\ref{ex3leaves} is defined on tha same surface $L$ as Example~\ref{ex2leaves} and the endperiodic map is defined as a composition $f=S\circ T$ of two homeomorphiams  $S$ and $T$ which are similar to the homeomorphisms $S$ and $T$ of Example~\ref{ex2leaves}.  In both examples $S$ and $T$ behave the same way near the ends of $L$.

 In Example~\ref{ex2leaves}, the handle in a fundamental domain of an end was considered as a hole through the surface (see Figure~\ref{four1}) and the homeomorphism $S$ (respectively $T$) draged this hole across the core when it moved this handle from the lower right (respectively lowert left) to upper left (respectively upper right) end of the surface $L$. Thus the juntures are distorted  on both the back and front of the surface of Figure~\ref{four1}.
 
 In Example~\ref{ex3leaves}, the handle in a fundamental domain of an end is considered as a handle attaced to the front of the surface (see Figure~\ref{four2}) and the homeomorphism $S$ (respectively $T$) drags this handle across the core when it moves this handle from the lower right (respectively lowert left) to upper left (respectively upper right) end of the surface $L$. Thus the juntures on only the front of Figure~\ref{four2} are distorted.  The result of this is that  in Example~\ref{ex3leaves}, both $\Lambda_{+}$ and $\Lambda_{-}$ have three leaves and the example has a pair of complementary principal regions with three arms.
 
 The positive lamination $\Lambda_{+}$ consists of three leaves $\lambda_{1},\lambda_{2},\lambda_{3}$ which are the border leaves of a positive principal region $P$.  The leaves $\lambda_{2}$ and $\lambda_{3}$ are isolated. The leaf $\lambda_{1}$ is semi-isolated, being approached from below by an arm of $P$ bordered by an end of $\lambda_{2}$ and an end of $\lambda_{3}$.  In the traintrack of Figure~\ref{four2} this arm is represented by the solid black curve that approaches $\lambda_{1}$ from below returning infinitely often to the core.  The leaf $\lambda_{1}$ is drawn in Figure~\ref{four2} as a horizontal line bordering the principal region $P$ from below. One end of the leaf $\lambda_{1}$ and one end of $\lambda_{2}$ border an escaping arm of the principal region $P$ as does the other end of the leaf $\lambda_{1}$ and one end of the leaf $\lambda_{3}$. The six vertices of the nucleus of the principal region $P$ are indicated by black dots. The three leaves of $\Lambda_{-}$ and the negative prncipal region are similarly drawn in Figure~\ref{four2}.

\end{example} 

\vfill\eject

\section{Another example defined without using Dehn twists}

\begin{figure}[t]
\begin{center}
\begin{picture}(300,290)(20,-270)
\rotatebox{270}{\scalebox{.9}{\includegraphics[width=300pt]{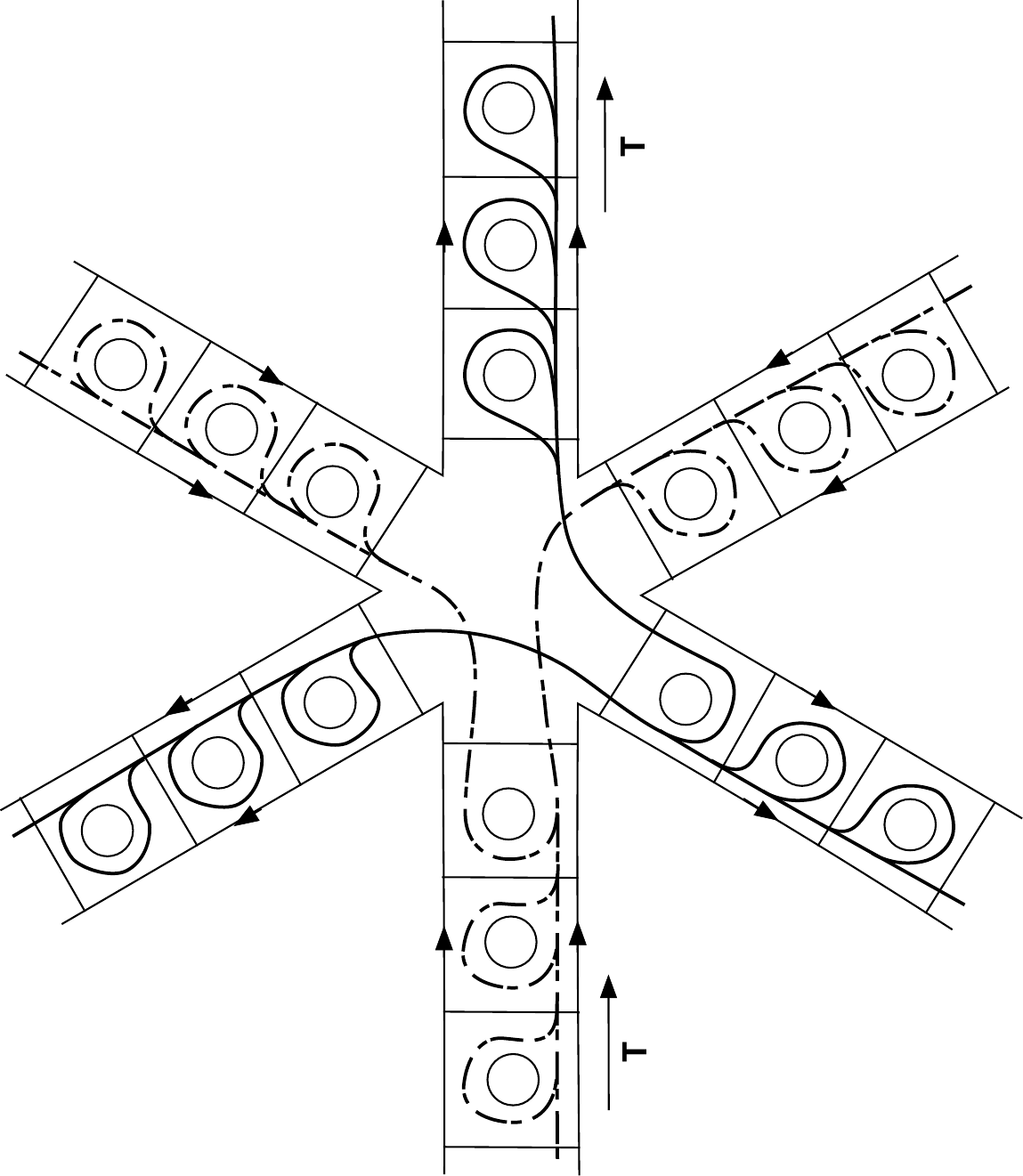}}}
\end{picture}
\caption{Traintracks for $\Lambda_{\pm}$ in Example~\ref{exsix}}\label{six}
\end{center}
\end{figure}

\begin{example}\label{exsix}

Let $L$ be the surface of Figure~\ref{six}.  Let $T$ be the homeomorphism that moves each circle boundary (and in fact the fundamental domains) to the right by one in the horizontal strips. The rightmost circle in the left strip is moved across the core to the leftmost circle in the right strip. The homeomorphism $T$ is to leave the other four strips pointwise fixed. Let $R_{3}$ be the rotation of the whole figure clockwise through $2\pi/3$ radians. Define $f = (R_{3}\circ T)^{3}$. There will be no fixed points on any of the boundary components. There will be three leaves with an escaping end in each of $\Lambda_{\pm}$ and these leaves will be the semi-isolated leaves in each of $\Lambda_{\pm}$.  Since each of the semi-isolated leaves accumulates on itself, both $\Lambda_{\pm}$ are transversally a Cantor set and there are uncountably many leaves in both $\Lambda_{\pm}$. Figure~\ref{six} illustrates  traintracks for $\Lambda_{\pm}$.

\end{example}

\begin{rem}

The method of Example~\ref{exsix} can be used to create examples of endperiodic automorphisms on planar surfaces with $2n$ nonsimple ends for any odd integer $n\ge 3$.

\end{rem}

\vfill\eject

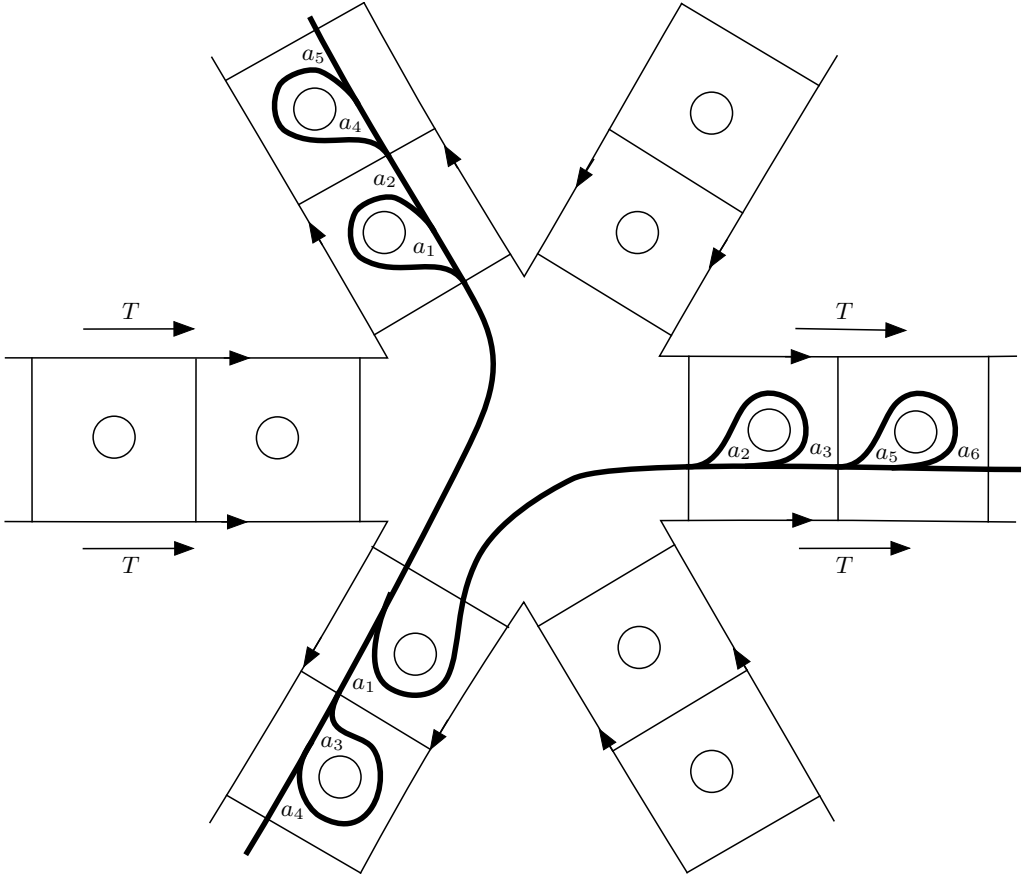
\begin{figure}[t]
\begin{center}
\begin{picture}(300,330)(45,-330)
\rotatebox{270}{\scalebox{1.1}{\includegraphics[width=300pt]{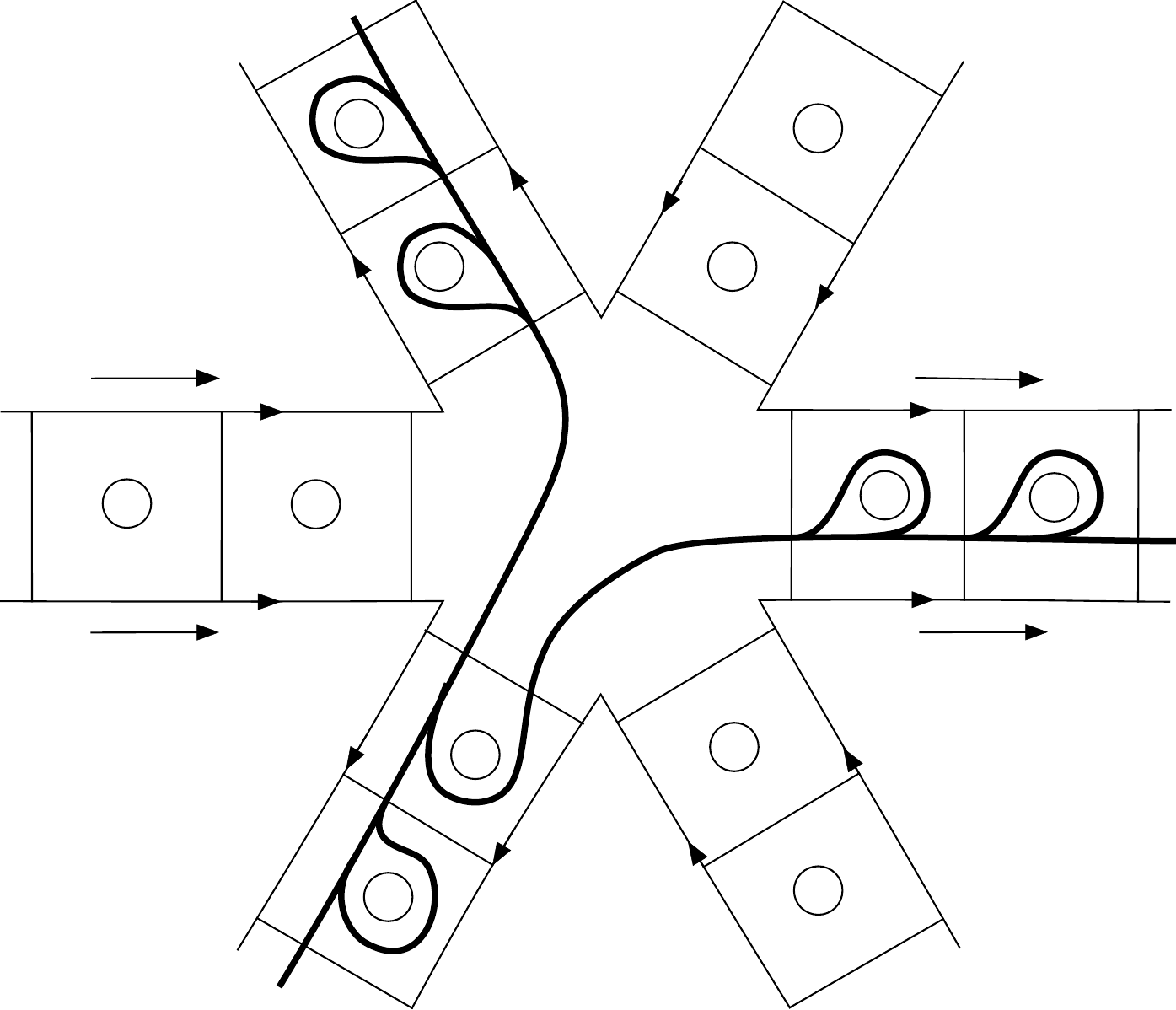}}}

\put (-253,-260){\small$a_{1}$}
\put (-245,-69){\small$a_{2}$}
\put (-80,-171){\small$a_{3}$}
\put (-280,-308){\small$a_{4}$}
\put (-272,-22){\small$a_{5}$}
\put (-24,-172){\small$a_{6}$}

\put (-230,-95){\small$a_{1}$}
\put (-111,-172){\small$a_{2}$}
\put (-265,-282){\small$a_{3}$}
\put (-258,-48){\small$a_{4}$}
\put (-55,-173){\small$a_{5}$}

\put (-340,-216){\small$T$}
\put (-340,-120){\small$T$}
\put (-70,-216){\small$T$}
\put (-70,-120){\small$T$}
\end{picture}
\caption{Traintracks for $\Lambda_{+}$ in Example~\ref{exsix}}\label{sixPlus}
\end{center}
\end{figure}

\begin{figure}[t]
\begin{center}
\begin{picture}(300,360)(55,-360)
\rotatebox{270}{\scalebox{1.2}{\includegraphics[width=300pt]{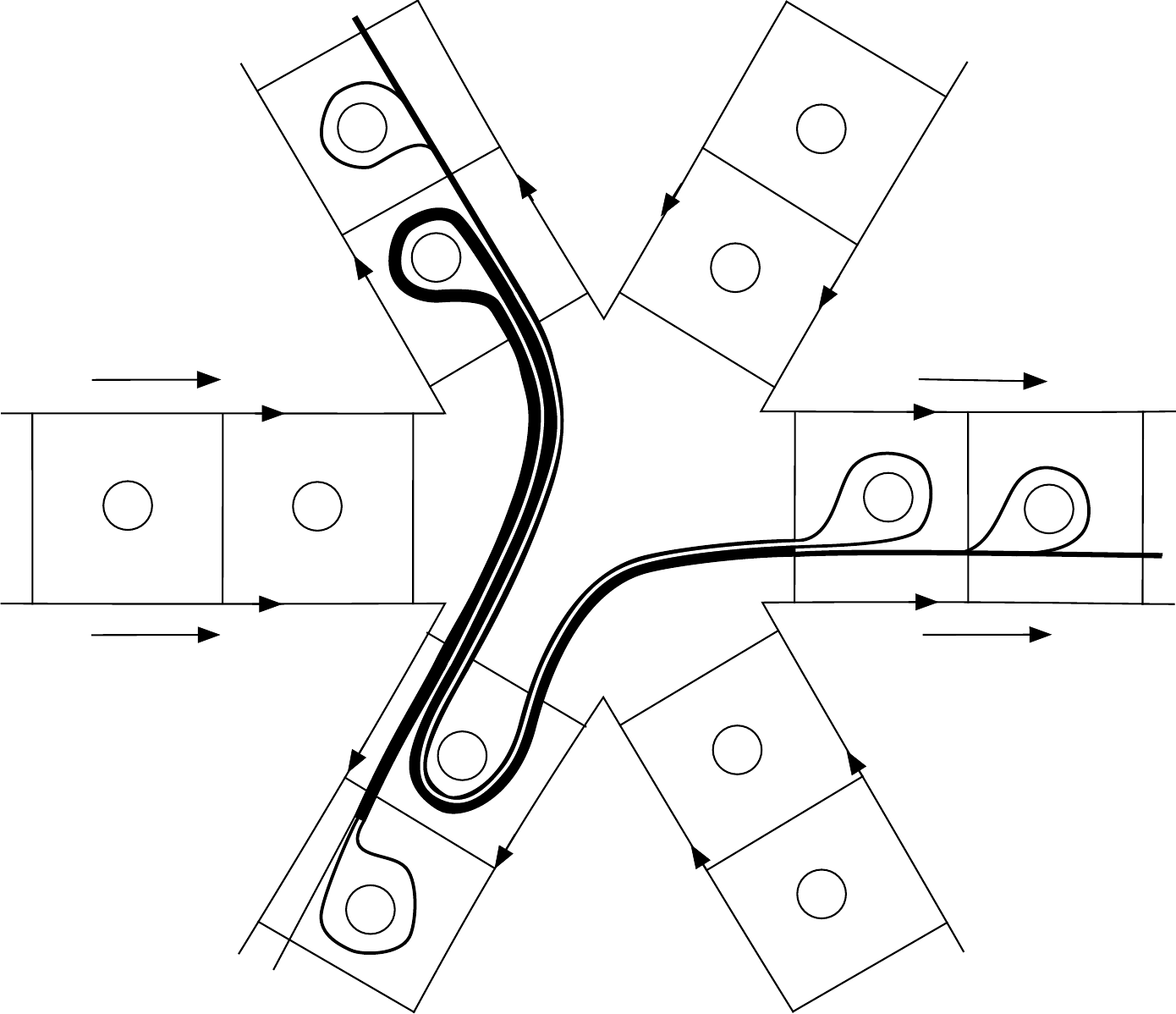}}}
\put (-280,-282){\small$a_{1}$}
\put (-265,-71){\small$a_{2}$}
\put (-90,-192){\small$a_{3}$}

\put (-245,-295){\small$1$}
\put (-285,-85){\small$x$}
\put (-252,-67){\small$x^{2}$}
\put (-88,-164){\small$x^{2}$}
\put (-268,-329){\small$x^{3}$}
\put (-310,-35){\small$x^{4}$}
\put (-35,-167){\small$x^{5}$}

\put (-255,-100){\small$a_{1}$}
\put (-123,-190){\small$a_{2}$}
\put (-295,-309){\small$a_{3}$}

\put (-370,-236){\small$T$}
\put (-370,-130){\small$T$}
\put (-70,-236){\small$T$}
\put (-70,-130){\small$T$}
\end{picture}
\caption{The traintrack $\TT_{2}$ in Example~\ref{exsix}}\label{sixB}
\end{center}
\end{figure}
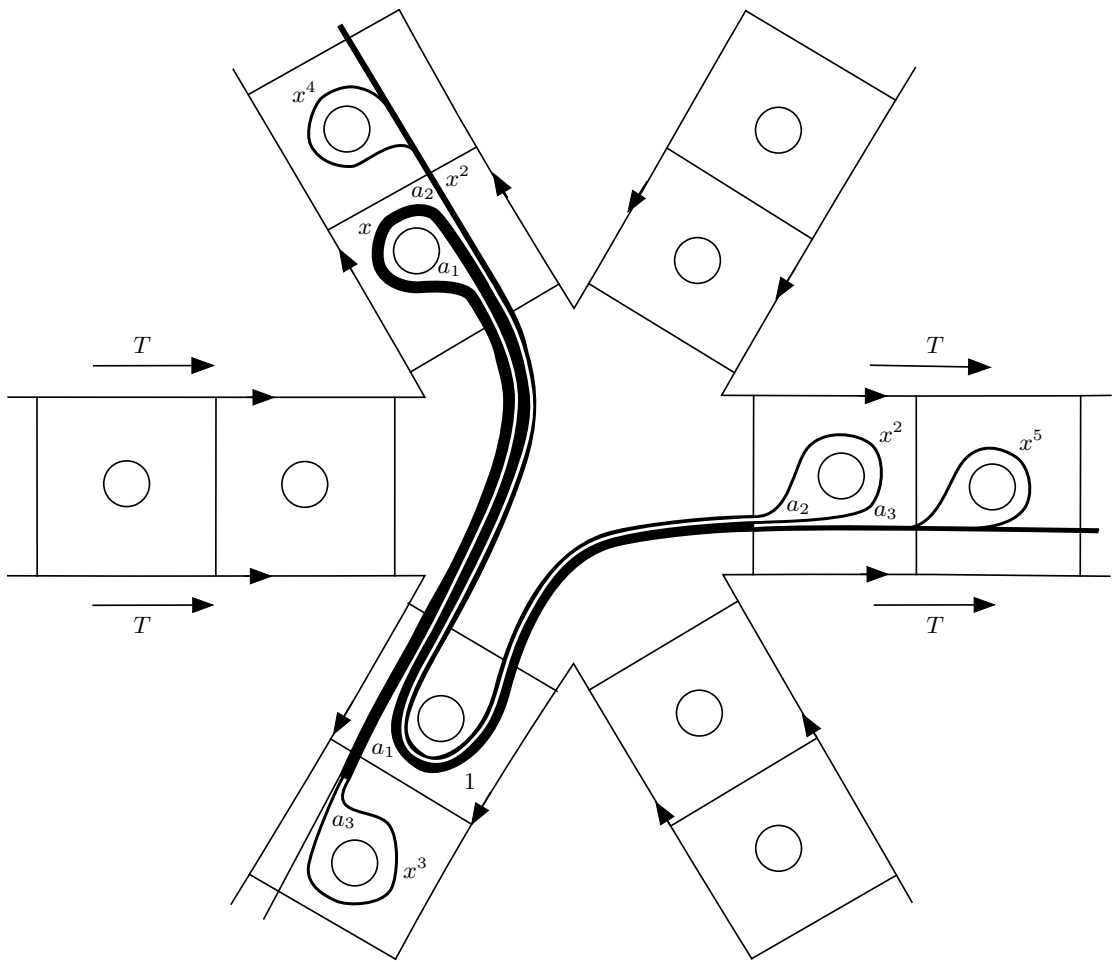

Denote the traintrack for $\Lambda_{+}$ given in Figure~\ref{six} by $\TT_{0}$. We have redrawn $\TT_{0}$ in Figure~\ref{sixPlus}. Let $g = R_{3}\circ T$. It is easier to work with $g$ than $f=g^{3}$. Let $\TT_{1} = g(\TT_{0})$, $\TT_{2} = g(\TT_{1})$, and in general $\TT_{n} = g(\TT_{n-1}) = g^{n}(\TT_{0})$, $n\ge 1$. The traintrack $\TT_{1}$ is obtained from $\TT_{0}$ by blowing air from $a_{1}$ to $a_{1}$ (see Figures~\ref{sixPlus} and~\ref{sixB}). Similarly $\TT_{2}$ is obtained from $\TT_{1}$ by blowing air from $a_{2}$ to $a_{2}$. Etc. In Figure~\ref{sixPlus} we have indicated how to blow air to create $\TT_{1},\TT_{2},\TT_{3},\TT_{4,}\TT_{5},\ldots$. From these it is easy to see the pattern. We draw $\TT_{2}$ in Figure~\ref{sixB} and indicate how to blow air to get $\TT_{3} = g^{3}(\TT_{0})$. 

Further, $|\Lambda_{+}| = \bigcap_{n=0}^{\infty}\TT_{n}$.

\begin{rem}

It really is clear that $\TT_{1}$ is obtained from $\TT_{0}$ by blowing air from $a_{1}$ to $a_{1}$. From this the pattern follows immediately.

\end{rem}

\begin{rem}

It is not at all directly obvious that $\TT_{3} = f(\TT_{0})$. That is why we look at $g$ rather than $f = g^{3}$.

\end{rem}

\vfill\eject

\begin{prop}

There exists a tansverse measure $\mu$ on $\Lambda_{+}$ so that if $\alpha$ is any arc transverse to $\Lambda_{+}$ then $\mu(g(\alpha)) = x\mu(\alpha)$ where $x=1/\tau$ with $\tau$ the golden number.

\end{prop}

\begin{proof}

We construct the measure directly. In Figure~\ref{sixB}, the symbols $1,x,x^{2},x^{3},x^{4},\ldots$ denotes the value of the measure of a transverse arc to $\Lambda_{+}$ at each of these points. Clearly, $\mu(g(\alpha)) = x\mu(\alpha)$. It remains to show that $\mu(\alpha + \beta) = \mu(\alpha) + \mu(\beta)$. For this it suffices that $1=x+x^{2}$ or $x^{2}+x-1 = 0$. Since $x$ must be positive,
$$x = \frac{-1+\sqrt{5}}{2} = 1/\tau.$$
\end{proof}

\begin{rem}

The lamination $\Lambda_{-}$ has a transverse measure with scale factor $\tau$.

\end{rem}

\vfill\eject

\section{An uncountable family of translations}

The mapping class group of a compact or a finite area surface is at most countable. In this section we give  uncountable many translations of the surface with two nonplanar ends. These translations  are not isotopic but  are in some sense the same.

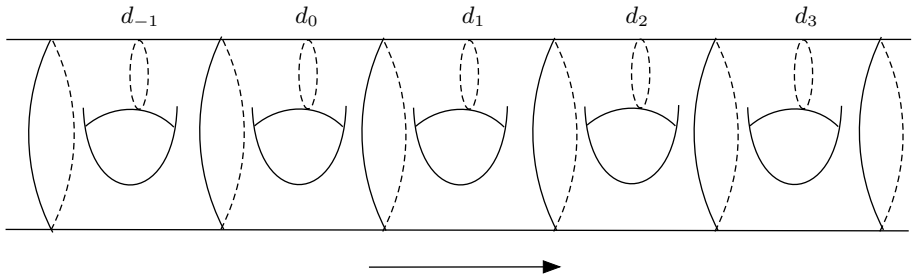
\begin{figure}[h]
\begin{center}
\begin{picture}(300,110)(20,-90)
\rotatebox{270}{\scalebox{.3}{\includegraphics[width=300pt]{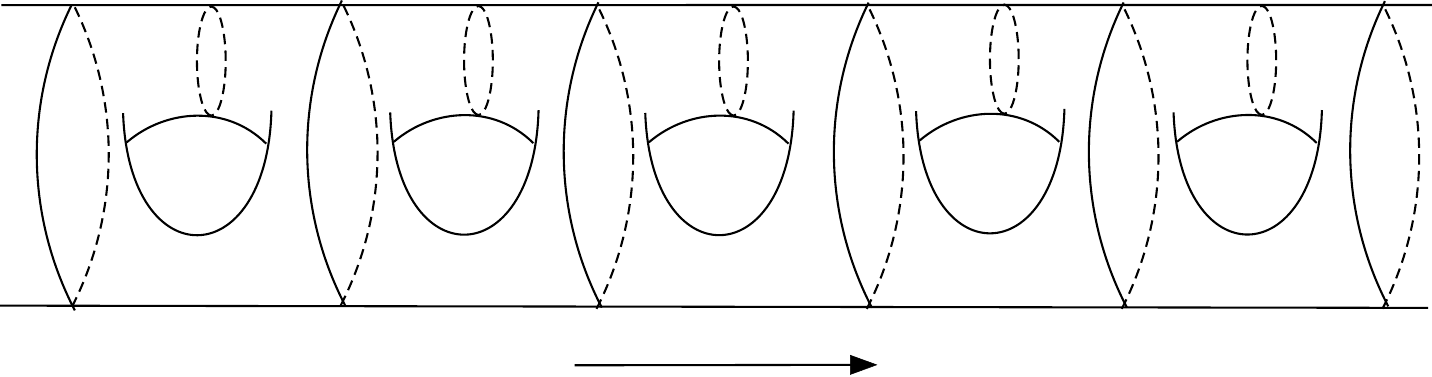}}}
\put (-299,5){\small$d_{-1}$}
\put (-233,5){\small$d_{0}$}
\put (-170,5){\small$d_{1}$}
\put (-108,5){\small$d_{2}$}
\put (-44,5){\small$d_{3}$}
\end{picture}
\caption{The surface $L$ admits uncountably many, nonisotopic translations}\label{transL}
\end{center}
\end{figure}

\begin{example}

Let $L$ be the surface of Figure~\ref{transL} with the indicated junctures and $d_{i}$, $-\infty<i<\infty$, be Dehn twists in the dotted circles. Let $T$ be any translation to the right with the indicated junctures. Let,
$$\iota = \{\ldots, i_{-1},i_{0},i_{1},i_{2},\ldots,i_{k},\ldots\}$$
be any sequence of zeroes and ones.
Define the homeomorphism $D_{\iota}:L\to L$ by 
$$D_{\iota} = \cdots d_{-1}^{i_{-1}}\circ d_{0}^{i_{0}}\circ d_{1}^{i_{1}}\circ d_{2}^{i_{2}}\circ\cdots\circ d_{k}^{i_{k}}\circ\cdots.$$
Then the family $D_{\iota}\circ T$ consists of uncountably many nonisotopic translations of $L$ with the same junctures.

\end{example}

\begin{rem}

If one allows different choices of the junctures, one gets even more isotopy classes of translations of $L$.

\end{rem}

\begin{rem}

Each of these translations can be realized as the monodramy of homeomorphic depth-one foliations of $S\times [0,1]$ where $S$ is the two holed torus.

\end{rem}

\section{An  uncountabe family of nonisotopic endperiodic automorphisms}

In this section we give another example of a surface $L$ having an uncountable family of nonisotopic endperiodic automorphisms. These examples have the property that for any two of them there exists a homeomorphism from $L$ to $L$ taking the positive and negative laminations $\Lambda_{\pm}$ for one to the positive and negative laminations for the other.

\begin{rem}

The technique can be used to give similar uncountable families of examples with much more complicated dynamics.

\end{rem}

\vfill\eject

\begin{figure}[htb]
\begin{center}
\begin{picture}(250,210)(20,-220)
\rotatebox{270}{\scalebox{.9}{\includegraphics[width=250pt]{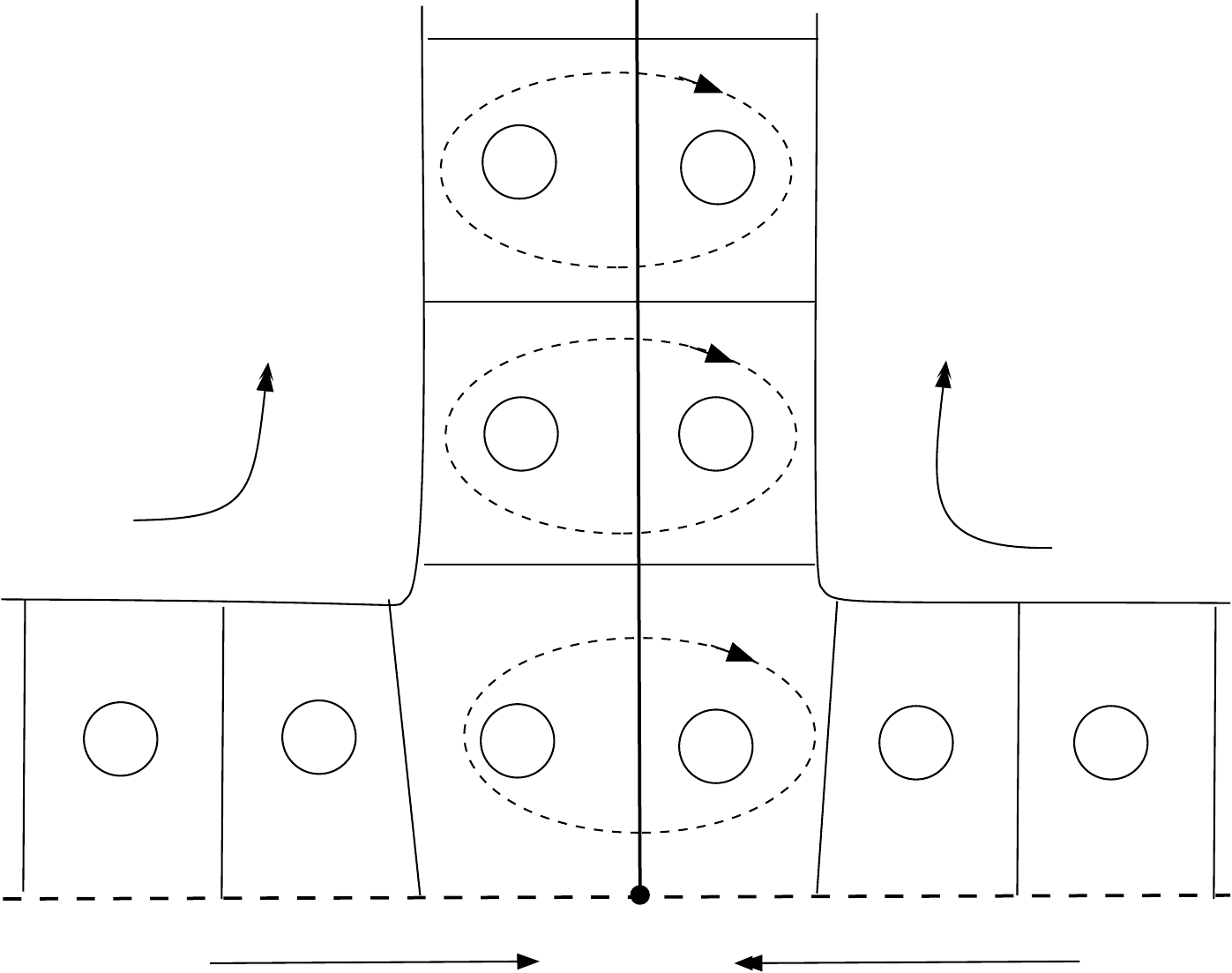}}}

\put (-115,-147){\small$d_{0}$}
\put (-113,-82){\small$d_{1}$}
\put (-115,-17){\small$d_{2}$}

\put (-155,-218){\small$\lambda_-$}
\put (-150,-143){\small$\lambda_{+}$}

\put (-289,-217){\small$J_{-2}$}
\put (-243,-217){\small$J_{-1}$}
\put (-193,-217){\small$J_0$}
\put (-10,-217){\small$J_{-2}$}
\put (-57,-217){\small$J_{-1}$}
\put (-102,-217){\small$J_0$}
\put (-203,-130){\small$J_1$}
\put (-203,-70){\small$J_2$}
\put (-203,-10){\small$J_3$}

\end{picture}
\caption{Top half of surface $L$ of Examples~\ref{ex2a}}\label{3PE}
\end{center}

\end{figure}
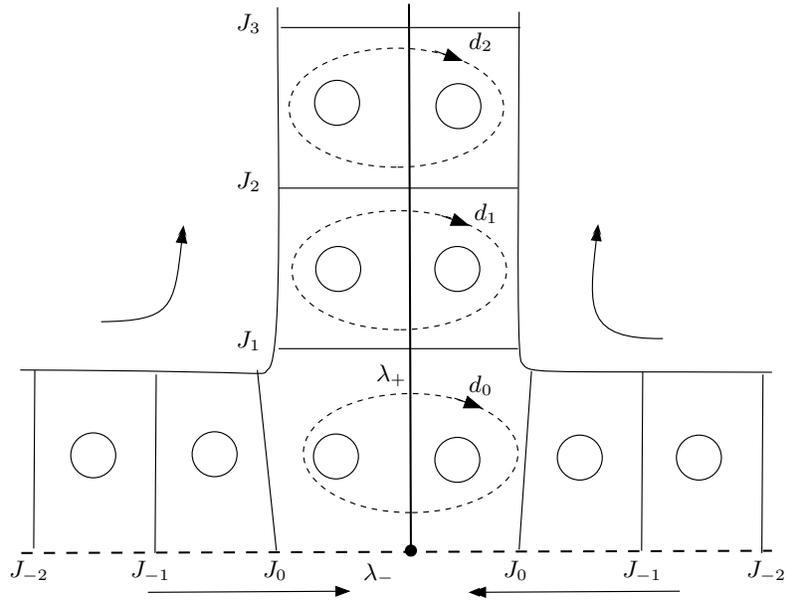

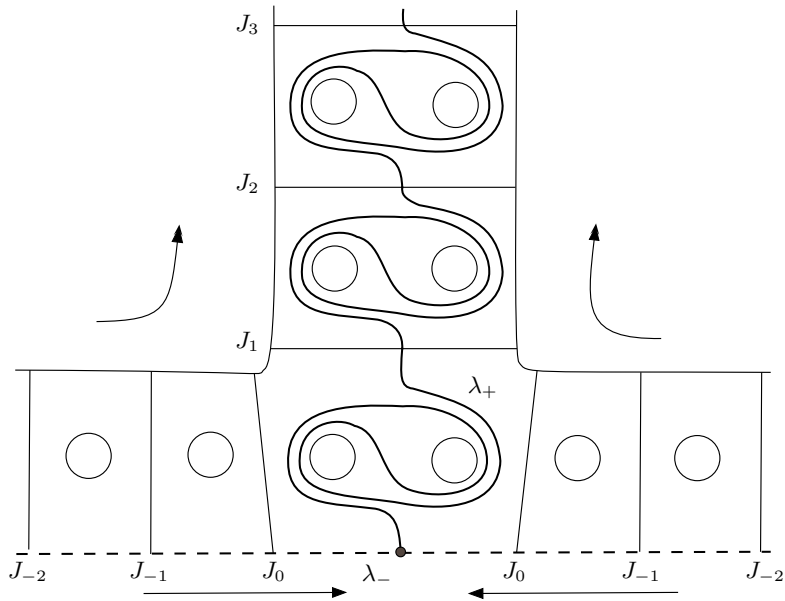
\begin{figure}[htb]
\begin{center}
\begin{picture}(250,210)(20,-220)
\rotatebox{270}{\scalebox{.9}{\includegraphics[width=250pt]{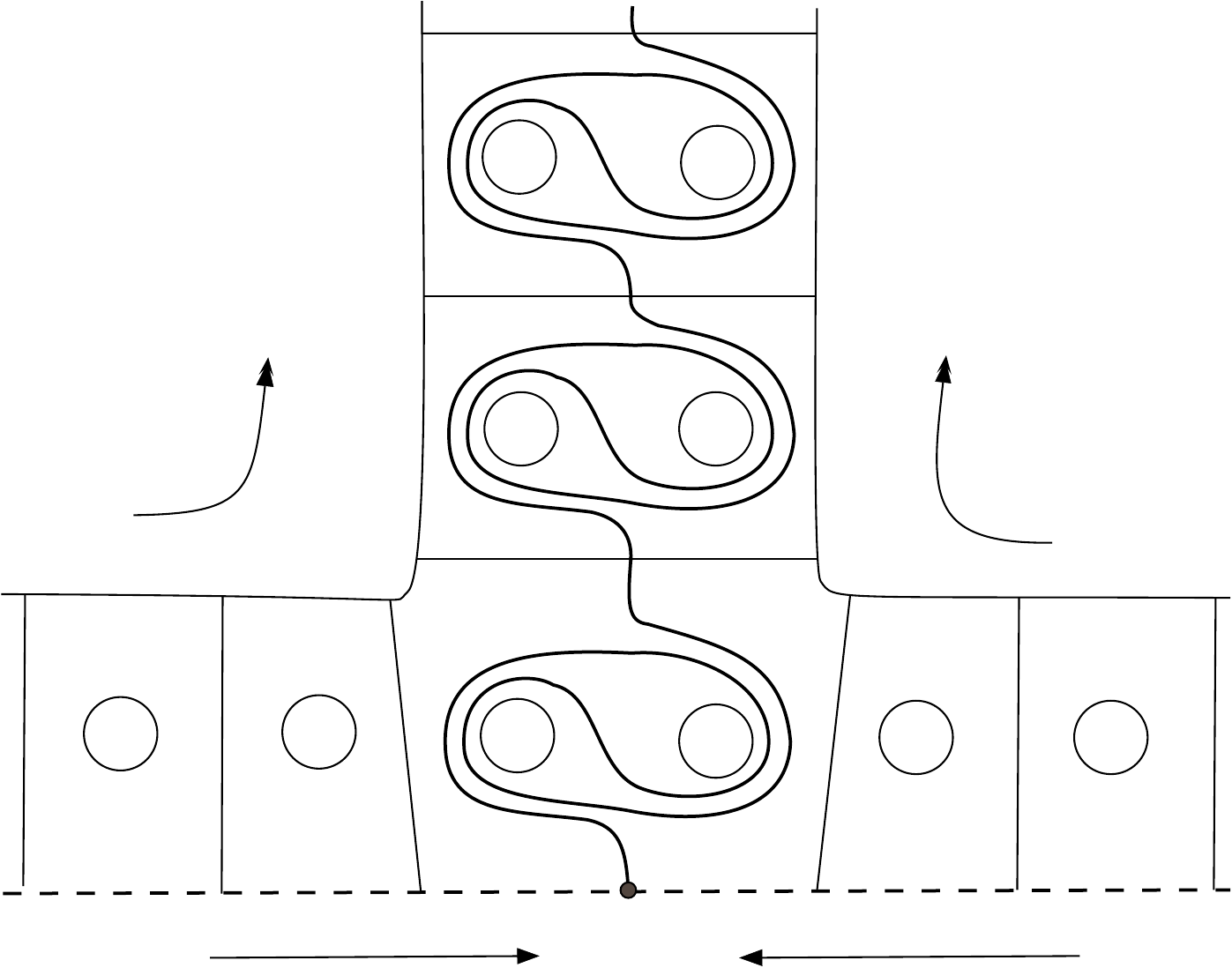}}}

\put (-155,-218){\small$\lambda_-$}
\put (-115,-147){\small$\lambda_{+}$}

\put (-289,-217){\small$J_{-2}$}
\put (-243,-217){\small$J_{-1}$}
\put (-193,-217){\small$J_0$}
\put (-10,-217){\small$J_{-2}$}
\put (-57,-217){\small$J_{-1}$}
\put (-102,-217){\small$J_0$}
\put (-203,-130){\small$J_1$}
\put (-203,-70){\small$J_2$}
\put (-203,-10){\small$J_3$}

\end{picture}
\caption{$\lambda_{+}$ when $\iota = \{1,0,0\ldots,0,\ldots\}$}\label{3PEA}
\end{center}

\end{figure}

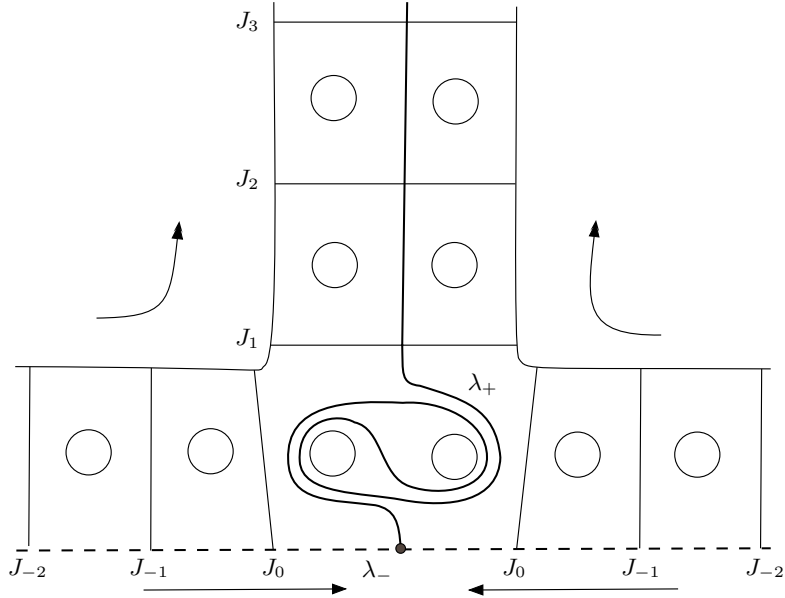
\begin{figure}[h]
\begin{center}
\begin{picture}(250,210)(20,-220)
\rotatebox{270}{\scalebox{.9}{\includegraphics[width=250pt]{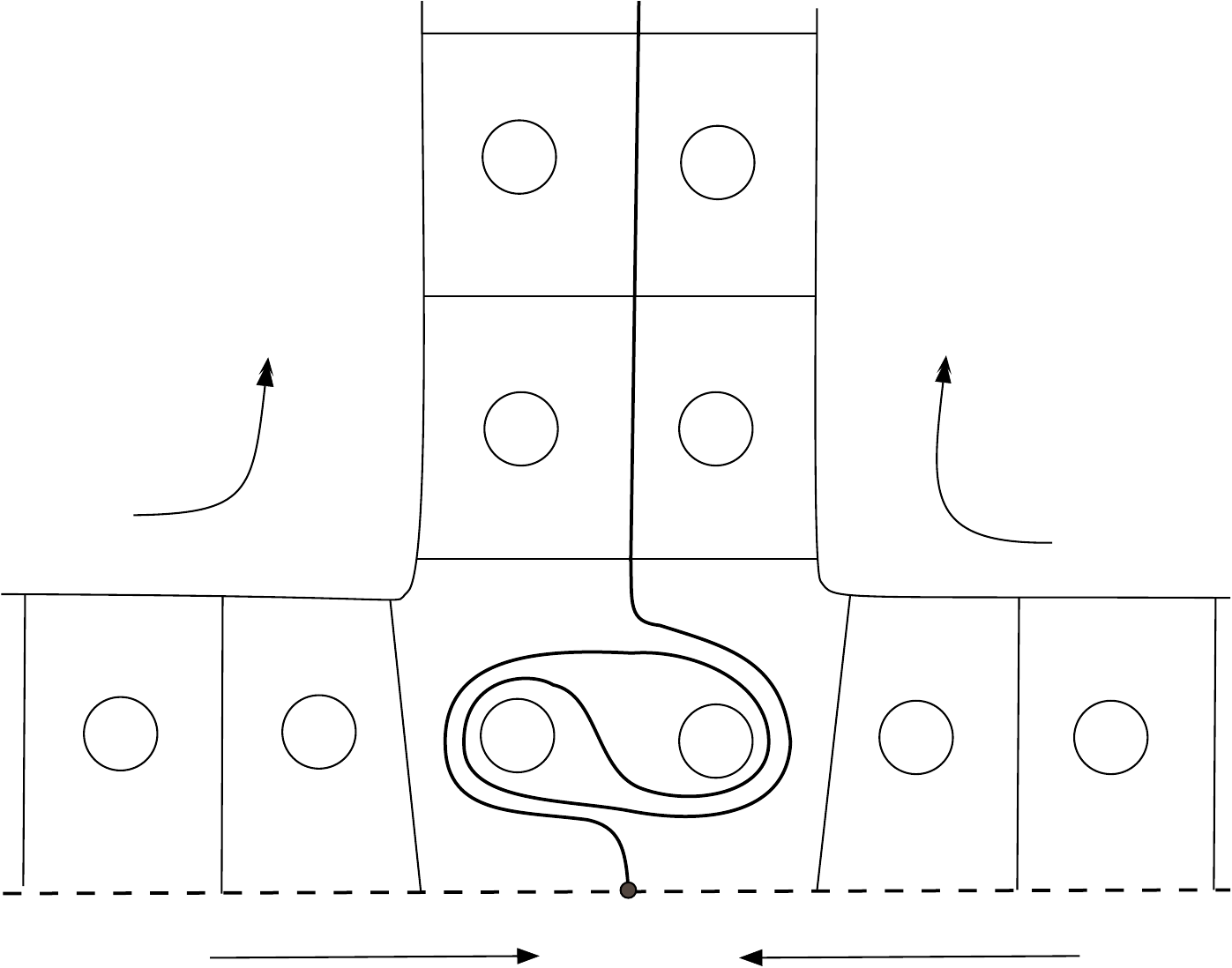}}}

\put (-155,-218){\small$\lambda_-$}
\put (-115,-147){\small$\lambda_{+}$}

\put (-289,-217){\small$J_{-2}$}
\put (-243,-217){\small$J_{-1}$}
\put (-193,-217){\small$J_0$}
\put (-10,-217){\small$J_{-2}$}
\put (-57,-217){\small$J_{-1}$}
\put (-102,-217){\small$J_0$}
\put (-203,-130){\small$J_1$}
\put (-203,-70){\small$J_2$}
\put (-203,-10){\small$J_3$}

\end{picture}
\caption{$\lambda_{+}$ when $\iota = \{1,-1,0,0\ldots,0,\ldots\}$}\label{3PEB}
\end{center}

\end{figure}

\begin{example}\label{ex2a}

Let $L$ be the planar surface of Figure~\ref{3PE} doubled along the bottom edge and $g:L\to L$  be an endperiodic automorphism with the two curves $\lambda_{+}$ and $\lambda_{-}$ of Figure~\ref{3PE} ae positive and negative laminations respectively.  We assume that   $g$   sends the junctures $J_i$ to $J_{i+1}$ for $-\infty < i < 0$ near the repelling ends and $1\le i < \infty$ near the attracting ends. There are two attracting ends and two repelling ends.  Let,
$$\iota = \{i_{0},i_{1},i_{2},\ldots,i_{k},\ldots\}$$
be any integer sequence.
Define the homeomorphism $D_{\iota}:L\to L$ by 
$$D_{\iota} =  \cdots\circ d_{k}^{i_{k}}\circ\cdots\circ d_{2}^{i_{2}}\circ d_{1}^{i_{1}}\circ d_{0}^{i_{0}}.$$
where $d_{i}$ is Dehn twist in the dotted curves indicated in Figure~\ref{3PE}. Then the family $D_{\iota}\circ g$ consists of uncountably many nonisotopic endperiodic automorphisms of $L$ such that,

\begin{enumerate}

\item  Each endperiodic automorphism has the $J_{i}$, $-\infty<i<\infty$, as junctures\upn{;}

\item $\Lambda_{-}$ and $\Lambda_{+}$ each consists of one isolated leaf.

\end{enumerate}

The leaf $\lambda_{-}$ is the bottom edge of Figures~\ref{3PE}, \ref{3PEA}, and~\ref{3PEB}. In Figures~\ref{3PEA} and~\ref{3PEB} we draw part of $\lambda_{+}$ for two different choices of $\iota = \{i_{0},i_{1},i_{2},\ldots,i_{k},\ldots\}$.

For all choices of $\iota = \{i_{0},i_{1},i_{2},\ldots,i_{k},\ldots\}$, the escaping set $\UU$ has four components each with border a curve of the first kind each of which gives a reducing geodesic. Reduction gives four planar strips with disks deleted approaching both ends and each having endperiodic automorphisms which are pure translations and one pseudo-anosov piece with four simple ends as in~\cite[Figure 17]{cc:HM}.

\end{example}

\section{An example with one repelling end and no attracting ends}

Our purpose in this section is to give an example of an endperiodic automorphism of a surface $L$ which has one repelling end $e$ and no attracting ends. By~\cite[Proposition~2.18]{cc:HM}, in such an example $L$ must have infinite endset. In fact by~\cite[Lemma~13.24]{cc:HM} and~\cite[Proposition~13.26]{cc:HM}, in such an example, the endset of $L$ must be uncountable. In our example the endset of $L$ will consist of the one isolated point $e$ and a Cantor set $C$.

We don't need the following proposition but it is the idea behind Example~\ref{cantorset}.

\begin{prop}

There exists a homeomorphism of the Cantor set without periodic points.

\end{prop}

\begin{proof}

We take as the Cantor set $C = \displaystyle\prod_{n=2}^{\infty}\Z_{n}$ where $\Z_{n}$ denotes the integers mod $n$. Define $\nu:C\to C$ by $\nu(\{a_{n}\}_{n=2}^{\infty}) = \{a_{n}+1\}_{n=2}^{\infty}$. Then $\nu$ is the desired homeomorphism.
\end{proof}

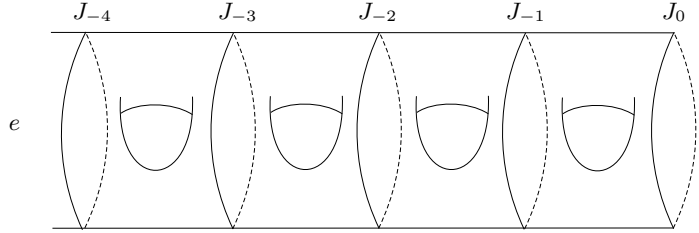
\begin{figure}[h]
\begin{center}
\begin{picture}(300,80)(-35,-70)
\rotatebox{270}{\scalebox{.25}{\includegraphics[width=300pt]{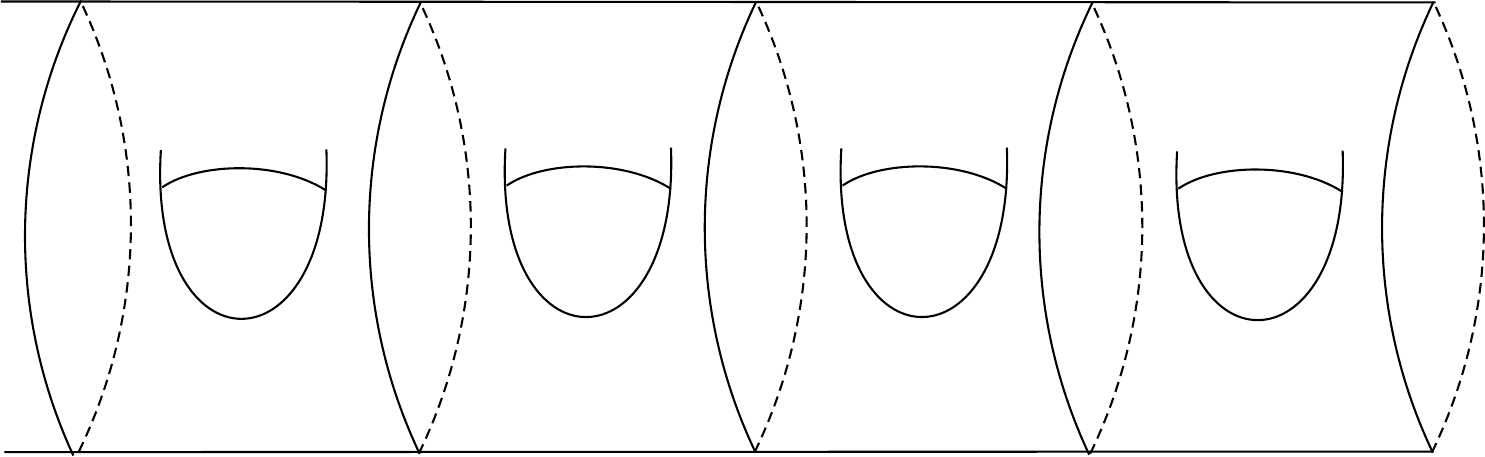}}}
\put (-260,-37){\small$e$}
\put (-236,5){\small$J_{-4}$}
\put (-181,5){\small$J_{-3}$}
\put (-128,5){\small$J_{-2}$}
\put (-72,5){\small$J_{-1}$}
\put (-13,5){\small$J_{0}$}
\end{picture}
\caption{The surface $L_{1}$ with end $e$ and boundary circle $J_{0}$}\label{L1}
\end{center}
\end{figure}

\begin{example}\label{cantorset}

We inductively construct the surface $L$ with endset consisting of one isolated end $e$  and a Cantor set $C$. To begin let $L_{1}$ be the surface with one circle boundary component and one nonplanar end $e$  in Figure~\ref{L1}. Let $S(n)$ denote the compact surface with $n+1$ boundary components and one handle. Let $L_{2}$ be $L_{1}$ with one copy of $S(2)$ attached along the one boundary component of $L_{1}$. let $L_{3}$ be $L_{2}$ with two copies of $S(3)$ attached along each of the two boundary components of $L_{2}$ and in general let $L_{n+1}$ be $L_{n}$ with $n!$ copies of $S(n+1)$ attached along the $n!$ boundary components of $L_{n}$. Let $L = \bigcup_{n=1}^{\infty}L_{n}$.

\end{example}

Note that the $L_{n}\sm\intr L_{n-1}$ consists of $(n-1)!$ copies of $S(n)$

\begin{lemma}

There exists a homeomorphism $g:L\to L$ such that,

\begin{enumerate}

\item $g$ fixes $L_{1}$ pointwise\upn{;}

\item $g(L_{n}) = L_{n}$, $n\ge 1$\upn{;}

\item If $S$ is a component of $L_{n}\sm\intr L_{n-1}$, then $g$ cyclically permutes the $n$ components of $\bd S\cap\bd L_{n}$. 

\end{enumerate}

\end{lemma}

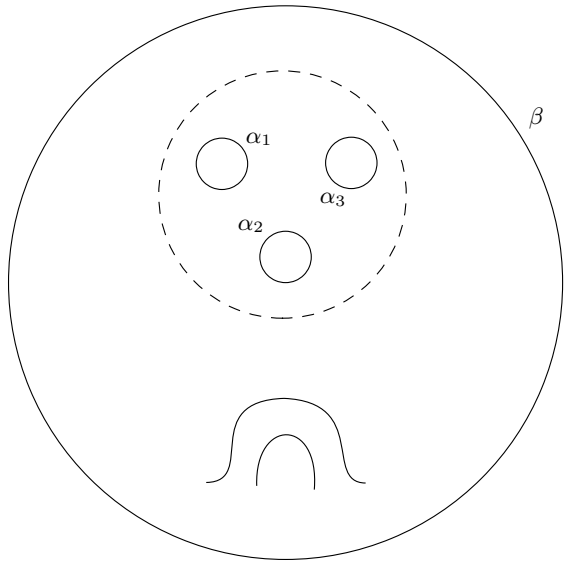
\begin{figure}[t]
\begin{center}
\begin{picture}(300,210)(-35,-210)
\rotatebox{270}{\scalebox{.7}{\includegraphics[width=300pt]{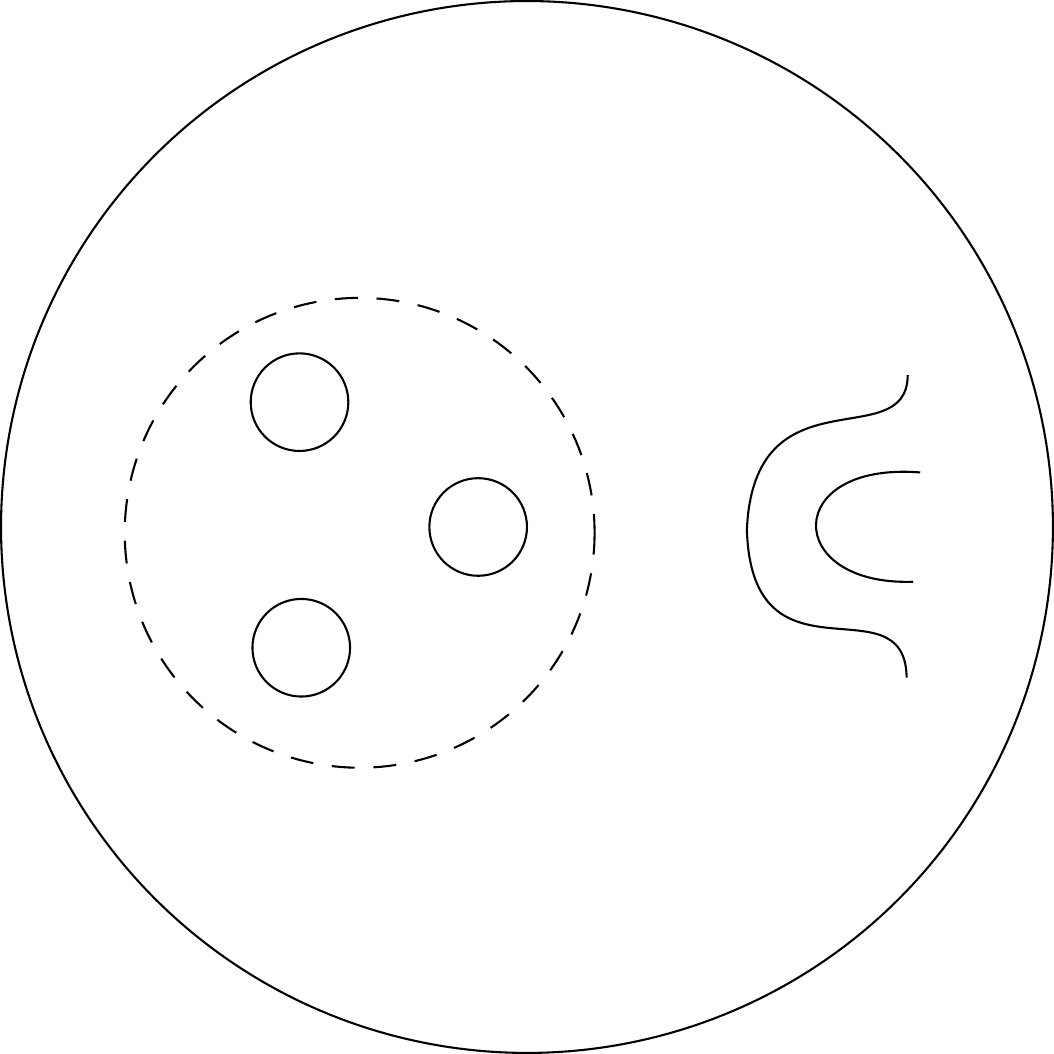}}}
\put (-120,-52){\small$\alpha_{1}$}
\put (-123,-85){\small$\alpha_{2}$}
\put (-92,-75){\small$\alpha_{3}$}
\put (-13,-45){\small$\beta$}
\end{picture}
\caption{The surface $S_{3}$ with 4 boundary circles  $\alpha_{1},\alpha_{2},\alpha_{3},\beta$ and one handle}\label{cyclic}
\end{center}
\end{figure}

\begin{proof}

We begin by setting $g = {\rm id}$ on $L_{1}$. If we have defined $g$ on $L_{n-1}$ and $S$ is a component of $L_{n}\sm\intr L_{n-1}$, then $S$ is homeomorphic to $S(n)$. Extend $g$ over $S$ by a homeomorphism supported in $\intr S$ that cyclically permutes the $n$ components of $\bd S\cap\bd L_{n}$. We illustrate such a homeomorphism for the case $n=3$ in Figure~\ref{cyclic}. The twist in the annulus (represented by the dotted circle in Figure~\ref{cyclic})  through $2\pi/3$ radians provides a homeomorphism from $S(3)$ to $S(3)$ which cyclically permutes three of the boundary components of $S(3)$. Inductively we have $g$ defined on $L$ satisfying the properties in the lemma. 
\end{proof}

\begin{lemma}

There exists a homeomorphism $h:L\to L$ such that,

\begin{enumerate}

\item $h$ fixes the endset of $L$ pointwise\upn{;}

\item $h(J_{i}) = J_{i+1}$, $i<0$ where the $J_{i}$ are the junctures in $L_{1}$ in Figure~\ref{L1}.

\end{enumerate}

\end{lemma}

\begin{proof}

For $n\ge 2$, choose a component $S_{n}$ of $L_{n}\sm\intr L_{n-1}$ homeomorphic to $S(n)$ such that $S_{n}\cap S_{n+1}$ is a circle. The homomorphism $h$ will be supported in $L_{1}\cup\bigcup_{n=2}^{\infty}S_{n}$. It is clear how to define $h$ to the left of the juncture $J_{-1}$ in Figure~\ref{L1}. Define the homeomorphism $h$ between $J_{-1}$ and $J_{0}$ to move the handle between $J_{-1}$ and $J_{0}$ into $S_{2}$. Define the homeomorphism $h$ on  $S_{2}$ to move the handle in $S_{2}$ into $S_{3}$ and inductively define the homeomorphism $h$ on  $S_{n}$ to move the handle in $S_{n}$ into $S_{n+1}$. The homeomorphism $h$ satisfies the properties of the lemma.
\end{proof}

The following proposition follows immediately from the two previous lemmas.

\begin{prop}

The endperiodic automorphism $f = h\circ g:L\to L$ has $e$ as repelling end and no other periodic end.

\end{prop}

\bibliographystyle{amsplain}

\begin{thebibliography}{10}


\bibitem{cc:HM}
J.~Cantwell and L.~Conlon, 
Endperiodic automorphisms of surfaces and foliations, 
arXiv:1006.4525v5.

\bibitem{fe:thesis}
S.~Fenley, \emph{Depth one foliations in hyperbolic $3$--manifolds}, {\it
  Thesis, Princeton University, 1990}.

\bibitem{fe:endp} 
S.~Fenley, \emph{Endperiodic surface
homeomorphisms and $3$--manifolds}, Math.  Z. 224 (1997) 1--24

\end{thebibliography}

\end{document}